\documentclass[12pt]{article}

\usepackage{graphicx}
\usepackage{amsmath}
\usepackage{amssymb}
\usepackage{amsfonts}
\usepackage{algorithm}
\usepackage[noend]{algpseudocode}
\title{A Treecode Algorithm for 3D Stokeslets and Stresslets}
\author{Lei Wang\thanks{Department of Mathematical Sciences, University of Wisconsin-Milwaukee, Milwaukee, WI 53211, USA wang256@uwm.edu} \and Svetlana Tlupova\thanks{Department of Mathematics, Farmingdale State College, SUNY, Farmingdale, NY 11735, USA tlupovs@farmingdale.edu} \and Robert Krasny\thanks{Department of Mathematics, University of Michigan, Ann Arbor, MI 48109, USA krasny@umich.edu}}

\begin{document}
\maketitle




%
%
%

\begin{abstract}
The Stokeslet and stresslet kernels are commonly used in
boundary element simulations and singularity methods for slow viscous flow.
Evaluating the velocity induced by a collection of Stokeslets and stresslets
by direct summation requires $O(N^2)$ operations,
where $N$ is the system size.
The present work develops a treecode algorithm for 3D Stokeslets and stresslets
that reduces the cost to $O(N\log N)$.
The particles are divided into a hierarchy of clusters,
and well-separated particle-cluster interactions are computed by a far-field Cartesian Taylor approximation. 
The terms in the approximation are contracted to promote efficient computation.
Serial and parallel results display the performance of the treecode for several test cases.
In particular the method has relatively simple structure
and
low memory usage,
and
this enhances parallel efficiency for large systems.
\end{abstract}

{\bf Keywords:} Stokeslet, stresslet, fast summation, treecode, Taylor approximation.




\section{Introduction}\label{intro}

The slow steady flow of an incompressible viscous fluid is governed by the Stokes equations,
\begin{subequations}
\begin{align}
\nabla^2 {\bf u} - \nabla p &= 0, \\[0.125em]
\nabla \cdot {\bf u} &= 0,
\end{align}
\end{subequations}
where $\bf u$ is  the fluid velocity, $p$ is the pressure, and the viscosity is taken to be unity.
Many applications in fluid dynamics are modeled as particle interactions in Stokes flow, 
including for example
particle-laden fluid jets~\cite{PignatelNicolasGuazzelliSaintillan}, 
vibrations in microfluidic crystals~\cite{BeatusTlustyBarZiv}, 
cilia- and flagella-driven flows~\cite{Drescher2010,Smith2009}, 
free-surface flows of liquid drops~\cite{NitscheParthasarathi}, 
and
vesicle flows~\cite{VeerapaneniRahimianBirosZorin},
among others.
The Stokeslet and stresslet kernels are fundamental solutions of the Stokes equations
given in 3D (up to a numerical prefactor) by
\begin{subequations}
\begin{align}
\label{Stokeslet}
S_{ij}(\bf{x,y}) &= \frac{\delta_{ij}}{|{\bf x} - {\bf y}|} + \frac{(x_i - y_i)(x_j - y_j)}{|{\bf x} - {\bf y}|^3}, \\[0.125em]
\label{stresslet}
T_{ijl} (\bf{x,y}) &= \frac{(x_i - y_i)(x_j - y_j)(x_l - y_l)}{|{\bf x} - {\bf y}|^5},
\end{align}
\end{subequations}
where $\delta_{ij}$ is the Kronecker delta,
${\bf x} = (x_1,x_2,x_3), {\bf y} = (y_1,y_2,y_3)$,
and
indices $i,j,l = 1:3$ represent Cartesian coordinates.
The Stokeslet and stresslet kernels are commonly used in boundary element simulations
and
singularity methods for slow viscous flow~\cite{Pozrikidis1992}.

The $i$th component of the velocity induced by a set of Stokeslets and stresslets is
\begin{align}
u_i({\bf x}^m) &= 
\sum_{{n = 1}\atop{n \ne m}}^N S_{ij}({\bf x}^m,{\bf x}^n)f_j^n + 
\sum_{{n = 1}\atop{n \ne m}}^N T_{ijl}({\bf x}^m,{\bf x}^n) h_j^n \nu_l^n, \quad i = 1,2,3,
\label{StoStre}
\end{align}
where 
${\bf x}^m$ is a target particle,
${\bf x}^n$ is a source particle,
$f_j^n$ is a force weight,
$h_j^n$ is a dipole weight, 
and $\nu_l^n$ are the components of a unit normal vector to a surface.
Note that the Stokeslet term has an implicit sum over $j=1:3$
and
the stresslet term has an implicit sum over $j,l=1:3$;
for clarity in some places below these sums will be written out explicitly.
Equation~(\ref{StoStre}) is written for the case in which the targets and sources coincide,
but it is straightforward to handle problems where they are disjoint.

Evaluating the velocity~(\ref{StoStre}) for $m = 1:N$ by direct summation 
requires $O(N^2)$ operations, 
which is prohibitively expensive when $N$ is large. 
The same issue arises for interacting point masses, point charges, and point vortices,
and
many fast summation methods have been developed to reduce the cost,
including 
particle-mesh methods~\cite{Darden1995,Hockney1988},
the Fast Multipole Method (FMM)~\cite{Greengard1987},
and
treecodes~\cite{BarnesHut}.
These methods reduce the operation count to $O(N\log N)$ or $O(N)$ in principle,
while introducing approximations.
The FMM and treecode use multipole expansions of particle clusters
(near-field and far-field for the FMM,
but only far-field for the treecode),
while particle-mesh methods interpolate the particle strengths to a grid
where often the fast Fourier transform (FFT) is used to compute the sum.

A number of these fast summation methods have been developed in the context of Stokes flow
including 
a particle-mesh Ewald technique~\cite{Saintillan2005,Sierou2001}
and
a pre-corrected FFT method~\cite{Wang2006}.
Several extensions of the FMM have also been developed for 
Stokes flow~\cite{CoronaGreengardRachhVeerapaneni2017,Gimbutas2015,Wang2007,Ying2004}.
In one implementation~\cite{Tornberg2008},
the Stokeslet and stresslet sums are decomposed into several terms which are
computed by the FMM for Coulomb interactions~\cite{Fu2000}. 
The kernel-independent FMM~\cite{MalhotraBiros2015,MalhotraBiros2016} 
has been applied to simulate swimming microorganisms~\cite{Rostami2016}
using regularized Stokeslets~\cite{Cortez2001,cortez-fauci-medovikov-05}.
Recently the Spectral Ewald (SE) method was developed for Stokes flow using 
Gaussian spreading functions~\cite{Klinteberg2016,Klinteberg2017}.

These developments significantly improve the capability of fast summation methods for Stokes flow,
but it is important to investigate different approaches
and
understand their properties,
especially as new types of many-core computing platforms become available
with new challenges to the goal of maintaining parallel efficiency.
In particular,
as the relative cost of memory access rises in comparison with arithmetic operations,
it is worthwhile to investigate fast summation methods with different memory requirements
and
communication patterns.

In this context 
the present work contributes a treecode algorithm
that reduces the cost of the Stokeslet and stresslet sums~(\ref{StoStre}) to $O(N\log N)$.
In a treecode,
the particles are divided into a hierarchy of clusters,
and well-separated particle-cluster interactions are computed by a far-field approximation,
while nearby interactions are computed directly~\cite{BarnesHut}.
The original treecode used monopole approximations,
but later work starting with the FMM showed the advantage of using
higher-order multipole approximations~\cite{Greengard1987}.
Here we retain the structure of the treecode,
but we employ higher-order Cartesian Taylor series
for the far-field
approximation~\cite{Draghicescu1995,Duan2001,
Lindsay2001,LiJohnstonKrasny2009,Wang2011,ambrose-siegel-tlupova-13}.
We derive novel expressions for the Taylor coefficients of the Stokeslet and stresslet kernels
in terms of the Taylor coefficients of the Coulomb potential,
and
the far-field approximation is contracted for efficient evaluation,
following an approach used for direct summation in the FMMLIB3D code~\cite{FMMLIB3D,Gimbutas2015}.
A key feature of the proposed treecode is its relatively simple structure
and
low memory usage,
which together can enhance parallel efficiency.

The paper is organized as follows. 
Section~2 explains the particle-cluster interaction on which the treecode is based.
Section~3 derives expressions for the Taylor coefficients of the Stokeslet and stresslet kernels.
Section~4 presents an efficient method for computing the particle-cluster approximations.
Section~5 describes the treecode algorithm in detail. 
Section~6 presents serial and parallel numerical results showing the performance of the treecode
in terms of accuracy, efficiency, and memory usage.
A summary is given in section 7.


\section{Particle-Cluster Interactions}
\label{sec:1}

\begin{figure}[htb]
\begin{center}
\includegraphics[width=0.65\linewidth]{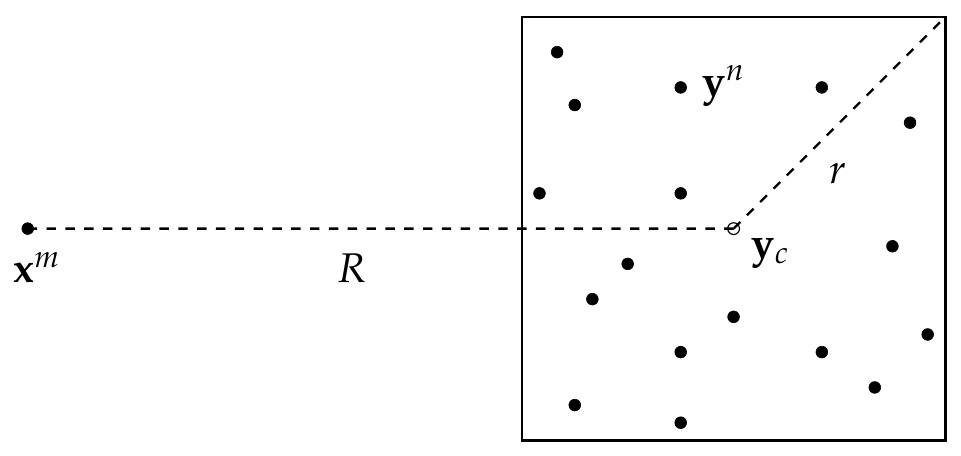}
\end{center}
\caption{
Particle-cluster interaction between a
target particle ${\bf x}^m$ 
and 
a cluster of source particles $C = \{ {\bf y}^n \}$,
with cluster center ${\bf y}_c$,
cluster radius $r$,
particle-cluster distance $R$.
}
\label{fig:particle_cluster}
\end{figure}

We start by expressing the velocity~(\ref{StoStre}) 
in terms of particle-cluster interactions.
Assume the particles have been divided into a set of clusters $\{ C \}$
(the procedure will be described below).
Then write the Stokeslet part of the velocity as a sum of particle-cluster interactions, 
\begin{subequations}
\begin{align}
u_i^{sto}({\bf x}^m) 
&= \sum_{n = 1}^{N}S_{ij}({\bf x}^m,{\bf x}^n)f_j^n \\
&= \sum_C\sum_{{\bf y}^n \in C}S_{ij}({\bf x}^m,{\bf y}^n)f_j^n =
\sum_C u_i^{sto}({\bf x}^m,C), 
\end{align}
\end{subequations}
where 
\begin{equation}
u_i^{sto}({\bf x}^m,C) = \sum_{{\bf y}^n \in C}S_{ij}({\bf x}^m,{\bf y}^n)f_j^n
\label{pc_definition}
\end{equation} 
is the interaction between a target particle ${\bf x}^m$
and 
a cluster of source particles $C = \{ {\bf y}^n \}$.
Figure~\ref{fig:particle_cluster} depicts the particle-cluster interaction,
showing also the 
cluster center ${\bf y}_c$,
cluster radius $r = \max_n |{\bf y}^n - {\bf y}_c|$,
and
particle-cluster distance $R = |{\bf x}^m - {\bf y}_c|$.
The stresslet part of the velocity is treated similarly.

If the particle ${\bf x}^m$ and cluster $C$ are not well-separated (the criterion is given later), 
then direct summation is used in~(\ref{pc_definition}).
If they are well-separated, 
then we use a far-field approximation
given by Taylor expanding the Stokeslet $S_{ij}({\bf x},{\bf y})$
about the cluster center ${\bf y} = {\bf y}_c$,
\begin{subequations}
\begin{align}
u_i^{sto}({\bf x}^m,C)
&= \sum_{{\bf y}^n \in C}S_{ij}({\bf x}^m,{\bf y}^n)f_j^n \\
&= \sum_{{\bf y}^n \in C}
\sum_{||{\bf k}|| = 0}^\infty
\frac{1}{{\bf k}!}D_{\bf y}^{\bf k }S_{ij}({\bf x}^m,{\bf y}_c)({\bf y}^n - {\bf y}_c)^{\bf k}f_j^n \\
&= \sum_{||{\bf k}|| = 0}^\infty
\frac{1}{{\bf k}!}D_{\bf y}^{\bf k}S_{ij}({\bf x}^m,{\bf y}_c)
\sum_{{\bf y}^n\in C}({\bf y}^n - {\bf y}_c)^{\bf k }f_j^n \\
&\approx \sum_{||{\bf k}|| = 0}^p a^{{\bf k}}_{ij}({\bf x}^m,{\bf y}_c)M_j^{\bf k}(C), 
\label{pc_interactionS}
\end{align}
\end{subequations}
where $p$ is the order of approximation, 
${\bf k} = (k_1,k_2,k_3)$ is an integer multi-index with all $k_i\geq 0$,
\begin{equation}
a^{\bf k}_{ij}({\bf x}^m,{\bf y}_c) =
\frac{1}{{\bf k}!}D_{\bf y}^{\bf k}S_{ij}({\bf x}^m,{\bf y}_c)
\label{stokeslet_a}
\end{equation}
is the ${\bf k}$th Taylor coefficient of the Stokeslet, 
and 
\begin{equation}
M_j^{{\bf k}}(C) = \sum_{{\bf y}^n \in C}({\bf y}^n - {\bf y}_c)^{\bf k }f_j^n
\end{equation}
is the ${\bf k}$th moment of cluster $C$. 
Note that the following conventions are used for multi-index notation,
$||{\bf k}|| = k_1 + k_2 + k_3,~
{\bf k}! = k_1!k_2!k_3!,~
D_{\bf y}^{\bf k} = D_{y_1}^{k_1}D_{y_2}^{k_2}D_{y_3}^{k_3},~
{\bf y}^{\bf k} = y_1^{k_1}y_2^{k_2}y_3^{k_3}$.

Note that the particle-cluster approximation~(\ref{pc_interactionS}) 
separates the target and source particles,
and
this promotes an efficient computation.
The key points are
(i) the Taylor coefficients $a_{ij}^{\bf k}$ depend on the target particle ${\bf x}^m$ 
and cluster center ${\bf y}_c$,
but not on the particles ${\bf y}^n$ in the cluster,
and
(ii) the cluster moments $M_j^{{\bf k}}$ depend only on the particles in the cluster $C$,
so they can be stored and re-used for different target particles ${\bf x}^m$. 
 
Similarly for the stresslet particle-cluster interaction,
\begin{equation}
u_i^{str}({\bf x}^m,C) =
\sum_{{\bf y}^n\in C}T_{ijl}({\bf x}^m,{\bf y}^n)h_j^n \nu_l^n \approx
\sum_{||{\bf k}|| = 0}^p\tilde{a}^{{\bf k}}_{ijl}({\bf x}^m,{\bf y}_c)\tilde{M}_{jl}^{{\bf k}}(C), 
\label{pc_interactionT}
\end{equation}
with stresslet Taylor coefficients,
\begin{equation}
\tilde{a}^{{\bf k}}_{ijl}({\bf x}^m,{\bf y}_c) = 
\frac{1}{{\bf k}!}  D_{{\bf y}}^{{\bf k}}T_{ijl}({\bf x}^m,{\bf y}_c),
\label{stresslet_a}
\end{equation}
and cluster moments,
\begin{equation}
\tilde{M}^{{\bf k}}_{jl}(C) = \sum_{{\bf y}^n\in C}({\bf y}^n - {\bf y}_c)^{{\bf k}} h_j^n \nu_l^n.
\end{equation}
Combining~(\ref{pc_interactionS}) and (\ref{pc_interactionT}),
when a particle ${\bf x}^m$ and cluster $C$ are well-separated, 
the far-field approximation for the induced velocity of the particle-cluster interaction is
computed from
\begin{equation}
u_i({\bf x}^m,C) \approx
\sum_{||{\bf k}|| = 0}^{p} a^{{\bf k}}_{ij}({\bf x}^m,{\bf y}_c)M_j^{{\bf k}}(C) +  
\sum_{||{\bf k}|| = 0}^{p}\tilde{a}^{{\bf k}}_{ijl}({\bf x}^m,{\bf y}_c)\tilde{M}_{jl}^{{\bf k}}(C).
\label{pc_interaction}
\end{equation}
As explained in the next two subsections,
the strategy for evaluating~(\ref{pc_interaction}) has two parts.
First we derive alternative expressions for the Stokeslet and stresslet Taylor coefficients,
$a^{{\bf k}}_{ij}, \tilde{a}^{{\bf k}}_{ijl}$,
and
then using those expressions we derive a method for efficient computation of~(\ref{pc_interaction}).
 

\section{Alternative expressions for Taylor coefficients}

First consider the Coulomb potential,
\begin{equation}
G({\bf x}, {\bf y}) = \frac{1}{|{\bf x} - {\bf y}|},
\label{GreenFunc}
\end{equation}
with Taylor coefficients
\begin{equation}
b^{{\bf k}}({\bf x}, {\bf y}) = \frac{1}{{\bf k}!}D_{{\bf y}}^{{\bf k}}G({\bf x},{\bf y}).
\end{equation}
It was shown in~ \cite{Duan2001,Lindsay2001} that these coefficients satisfy the recurrence relation
\begin{equation}
\label{recurrence_1}
||{\bf k}||\cdot |{\bf x} - {\bf y}|^2b^{{\bf k}} - (2 ||{\bf k}|| - 1)\sum_{i = 1}^{3}(x_i - y_i)b^{{\bf k} - {\bf e}_i} + (||{\bf k}|| - 1)\sum_{i = 1}^{3}b^{{\bf k} - 2{\bf e}_i} = 0, 
\end{equation}
with $b^{\bf{0}} = G({\bf x},{\bf y})$, $b^{{\bf k}} = 0$ when any $k_i < 0$, 
and ${\bf e}_i$ is the $i$th Cartesian unit vector. 
First the Stokeslet~(\ref{Stokeslet}) is written as 
\begin{equation}
S_{ij}({\bf x},{\bf y}) = 
\delta_{ij} G({\bf x},{\bf y}) + (x_{j}-y_{j})D_{y_i} G({\bf x},{\bf y}),
\end{equation}
and
then applying the Leibniz rule for differentiating a product, we have
\begin{subequations}
\begin{align}
\frac{1}{{\bf k}!} D_{\bf y}^{{\bf k}}&\left[(x_{j}-y_{j})D_{y_i}G\right] =
\frac{1}{{\bf k}!} D_{\bf y}^{{\bf k}-k_{j}{\bf e}_{j}} D_{y_{j}}^{k_{j}}\left[(x_{j}-y_{j}) D_{y_{i}} G\right] \\[4pt]
&= \frac{1}{{\bf k}!} D_{\bf y}^{{\bf k}-k_{j}{\bf e}_{j}} \left[(x_{j}-y_{j}) D_{y_{j}}^{k_{j}} D_{y_{i}} G - k_{j}D_{y_{j}}^{k_{j}-1} D_{y_{i}} G\right] \\[4pt]
&= (x_{j}-y_{j}) \frac{1}{{\bf k}!} D_{\bf y}^{{\bf k}+{\bf e}_{i}} G -
\frac{k_{j}}{{\bf k}!} D_{\bf y}^{{\bf k}+{\bf e}_{i}-{\bf e}_{j}} G \\[8pt]
&= (x_{j}-y_{j})(k_{i}+1)b^{{\bf k}+{\bf e}_{i}} - (k_{i} +1 -\delta_{ij}) b^{{\bf k}+{\bf e}_{i}-{\bf e}_{j}}.
\end{align}
\end{subequations}
This yields the following expression for $a^{{\bf k}}_{ij}$ in terms of the $b^{{\bf k}}$,
\begin{equation}
\label{recurrence_2}
a^{{\bf k}}_{ij} = 
\delta_{ij}b^{{\bf k}} + 
(x_j - y_j)(k_i + 1)b^{{\bf k} + {\bf e}_i} - 
(k_i +1 - \delta_{ij})b^{{\bf k} + {\bf e}_i-{\bf e}_j}.
\end{equation}
Similarly the stresslet~(\ref{stresslet}) is written as
\begin{equation}
T_{ijl}({\bf x},{\bf y}) = 
\frac{1}{3} \left[(x_l-y_l)D_{y_i}D_{y_j}G({\bf x},{\bf y}) + 
\delta_{ij}D_{y_l}G({\bf x},{\bf y})\right],
\end{equation}
so the Taylor coefficients are
\begin{equation}
\tilde{a}^{{\bf k}}_{ijl} =
\frac{1}{{\bf k}!} D_{\bf y}^{{\bf k}} T_{ijl} = \frac{1}{3} \left[\frac{1}{{\bf k}!} D_{\bf y}^{{\bf k}}  
\left[(x_{l}-y_{l})D_{y_i}D_{y_j}G\right] + \delta_{ij} \frac{1}{{\bf k}!} D_{\bf y}^{{\bf k}}D_{y_l}G\right].
\end{equation}
The first term on the right is
\begin{subequations}
\begin{align}
\frac{1}{{\bf k}!}D_{\bf y}^{{\bf k}}&\left[(x_{l} - y_{l}) D_{y_{i}}D_{y_{j}}G\right] \\
&= (x_{l}-y_{l}) \frac{1}{{\bf k}!} D_{\bf y}^{{\bf k}+{\bf e}_{i}+{\bf e}_{j}} G - 
\frac{k_{l}}{{\bf k}!} D_{\bf y}^{{\bf k}+{\bf e}_{i}+{\bf e}_{j}-{\bf e}_{l}} G \\[6pt]
&= (x_{l}-y_{l})(k_{i}+1)(k_{j} +1 +\delta_{ij}) b^{{\bf k}+{\bf e}_{i}+{\bf e}_{j}} \\[6pt]
& - (k_{i}+1-\delta_{il})(k_{j} +1+\delta_{ij}-\delta_{jl}) b^{{\bf k}+{\bf e}_{i}+{\bf e}_{j}-{\bf e}_{l}}.
\end{align}
\end{subequations}
This yields the following expression for $\tilde{a}^{{\bf k}}_{ijl}$ in terms of the $b^{{\bf k}}$,
\begin{align}
3\tilde{a}^{{\bf k}}_{ijl} &= (x_l - y_l)(k_i + 1)(k_j + 1 + \delta_{ij})b^{{\bf k} + {\bf e}_i + {\bf e}_j} \nonumber\\
& - (k_i + 1 - \delta_{il})(k_j + 1 + \delta_{ij} - \delta_{jl})b^{{\bf k} + {\bf e}_i + {\bf e}_j - {\bf e}_l} \nonumber\\
\label{recurrence_3}
& + \delta_{ij}(k_l + 1)b^{{\bf k} + {\bf e}_l}.
\end{align}
The Taylor coefficients $a^{{\bf k}}_{ij}, \tilde{a}^{{\bf k}}_{ijl}$ could be computed explicitly 
using (\ref{recurrence_2}) and (\ref{recurrence_3}),
however instead we will employ these relations implicitly to obtain a more efficient evaluation of the particle-cluster approximation~(\ref{pc_interaction}).
The details are explained in the next section.
   
 
\section{Efficient Computation of Particle-Cluster Approximations} 

In this section we explain how to rewrite various sums for improved computational efficiency,
first for the direct sum~(\ref{StoStre}) to illustrate the idea,
and
then for the particle-cluster approximation~(\ref{pc_interaction}).

\subsection{Direct sum}

The procedure for efficient evaluation of the direct sum~(\ref{StoStre}) 
was previously used in the FMMLIB3D code~\cite{FMMLIB3D,Gimbutas2015}.
It is repeated here to illustrate the idea,
and 
then it is applied to the particle-cluster approximation in the treecode.
The idea is to contract the sums and re-use quantities wherever possible.
First note that the Stokeslet part of~(\ref{StoStre}) can be expressed as
\begin{subequations}
\begin{align} 
u_i^{sto}({\bf x}^m)
&= \sum_{{n = 1}\atop{n \ne m}}^N\sum_{j=1}^3 S_{ij}({\bf x}^m,{\bf y}^n)f_{j}^{n} \\
&= \sum_{{n = 1}\atop{n \ne m}}^N
\left(
\sum_{j=1}^3 \frac{\delta_{ij}}{|{\bf x}^m - {\bf y}^n|}f_{j}^{n} + 
\sum_{j=1}^3 \frac{(x_i^m - y_i^n)(x_j^m - y_j^n)}{|{\bf x}^m - {\bf y}^n|^3}f_{j}^{n}\right) \\
&= \sum_{{n = 1}\atop{n \ne m}}^N
\left(
\frac{f_i^n}{|{\bf x}^m - {\bf y}^n|} + 
(x_i^m - y_i^n)s^{mn}\right),
\label{stokeslet_direct_sum}
\end{align}
\end{subequations}
where
$s^{mn} = |{\bf x}^m - {\bf y}^n|^{-3}\sum\limits_{j=1}^3 (x_j^m - y_j^n)f_j^n$
can be re-used for $i = 1:3$. 
Similarly, the stresslet part of~(\ref{StoStre}) can be expressed as
\begin{subequations}
\begin{align} 
u_i^{str}({\bf x}^m)
&= \sum_{{n = 1}\atop{n \ne m}}^N\sum_{j=1}^3\sum_{l=1}^3 T_{ijl}({\bf x}^m,{\bf y}^n) h_j^n \nu_l^n \\
&= \sum_{{n = 1}\atop{n \ne m}}^N\sum_{j=1}^3\sum_{l=1}^3 
\frac{(x_i^m - y_i^n)(x_j^m - y_j^n)(x_l^m - y_l^n)}{|{\bf x}^m - {\bf y}^n|^5}
h_j^n \nu_l^n \\
&= \sum_{{n = 1}\atop{n \ne m}}^N
(x_i^m - y_i^n)t^{mn},
\label{stresslet_direct_sum}
\end{align}
\end{subequations}
where
$t^{mn} = 
|{\bf x}^m - {\bf y}^n|^{-5}\sum\limits_{j=1}^3 (x_j^m - y_j^n)h_j^n\sum\limits_{l=1}^3 (x_l^m - y_l^n) \nu_l^n$
can be re-used for $i = 1:3$.
Using~(\ref{stokeslet_direct_sum}) and (\ref{stresslet_direct_sum}),
the operation count for direct summation is still $O(N^2)$,
but we observe empirically that the CPU run time is reduced
since this procedure avoids explicitly forming the tensors $S_{ij}, T_{ijl}$.


\subsection{Stokeslet particle-cluster interaction}

We take a similar approach in computing the particle-cluster approximations in the treecode;
the sums are contracted to reduce the number of operations,
terms are re-used wherever possible,
and
we avoid explicitly forming the Taylor coefficient tensors $a_{ij}^{\bf k}, \tilde{a}_{ijl}^{\bf k}$.
Consider the Stokeslet particle-cluster approximation~(\ref{pc_interaction});
using the expression for the Taylor coefficients~(\ref{recurrence_2}),
we have
\begin{subequations}
\begin{align}
u_i^{sto}({\bf x}^m,C) \approx
\sum_{||{\bf k}||=0}^p \sum_{j=1}^3 & a^{{\bf k}}_{ij}({\bf x}^m,{\bf y}_c) M_j^{{\bf k}}(C) \\
= \sum_{||{\bf k}||=0}^p \sum_{j=1}^3 & \Big[\delta_{ij}b^{{\bf k}} + (x_j - y_j)(k_i + 1)b^{{\bf k} + {\bf e}_i} \\
& - (k_i +1 - \delta_{ij})b^{{\bf k} + {\bf e}_i-{\bf e}_j}\Big] M_j^{{\bf k}}(C).
\end{align}
\end{subequations}
This simplifies to
\begin{equation}
\label{S_final}
u_i^{sto}({\bf x}^m,C) \approx\!\!
\sum_{||{\bf k}||=0}^p \Big[ 2 b^{{\bf k}} M_i^{{\bf k}}(C)
+ 
(k_i + 1) \Big[b^{{\bf k} + {\bf e}_i} \sigma^{{\bf k}}(C) - \sum_{j=1}^3 b^{{\bf k} + 
{\bf e}_i-{\bf e}_j} M_j^{{\bf k}}(C)\Big]\Big],
\end{equation}
where $\sigma^{{\bf k}}(C) = \sum\limits_{j=1}^3 (x_j - y_j) M_j^{{\bf k}}(C)$ can be re-used for $i=1:3$.


\subsection{Stresslet particle-cluster interaction}

Next consider the stresslet particle-cluster approximation~(\ref{pc_interaction});
using the expression for the Taylor coefficients~(\ref{recurrence_3}), 
we have
\begin{subequations}
\begin{align}
u_i^{str}&({\bf x}^m,C)
\approx \sum_{||{\bf k}||=0}^p \sum_{j=1}^3 \sum_{l=1}^3 
\tilde{a}^{{\bf k}}_{ijl}({\bf x}^m,{\bf y}_c)\tilde{M}_{jl}^{{\bf k}}(C) \\
&= \frac{1}{3}\sum_{||{\bf k}||=0}^p \sum_{j=1}^3 \sum_{l=1}^3 
\Big[(x_l - y_l)(k_i + 1)(k_j + 1 + \delta_{ij})b^{{\bf k} + {\bf e}_i + {\bf e}_j} \\
\label{T_1c}
& - (k_i + 1 - \delta_{il})(k_j + 1 + \delta_{ij} - \delta_{jl})b^{{\bf k} + {\bf e}_i + {\bf e}_j - {\bf e}_l} \\[0.5em]
& + \delta_{ij}(k_l + 1)b^{{\bf k} + {\bf e}_l} \Big] \tilde{M}_{jl}^{{\bf k}}(C).
\end{align}
\end{subequations}
To keep the next few intermediate formulas more concise,
in the remainder of this section we drop the arguments ${\bf x}^m, C$.
Then moving the sums over indices $j,l$ as far as possible to the right,
and
splitting~(\ref{T_1c}) into two terms, 
we have
\begin{subequations}
\begin{align}
u_i^{str} &\approx
\frac{1}{3}\sum_{||{\bf k}||=0}^p
\label{T_2a}
\bigg[ (k_i + 1) \sum_{j=1}^3 (k_j + 1 + \delta_{ij})b^{{\bf k} + {\bf e}_i + {\bf e}_j} 
\sum_{l=1}^3 (x_l - y_l)\tilde{M}_{jl}^{{\bf k}} \\
\label{T_2b}
& - (k_i + 1) \sum_{j=1}^3 \sum_{l=1}^3 (k_j + 1 + \delta_{ij} - \delta_{jl})b^{{\bf k} + {\bf e}_i + {\bf e}_j - {\bf e}_l}\tilde{M}_{jl}^{{\bf k}} \\
\label{T_2c}
& + \sum_{j=1}^3 \sum_{l=1}^3 \delta_{il} (k_j + 1 + \delta_{ij} - \delta_{jl})b^{{\bf k} + {\bf e}_i + {\bf e}_j - {\bf e}_l} \tilde{M}_{jl}^{{\bf k}} \\
\label{T_2d}
& + \sum_{l=1}^3 (k_l + 1)b^{{\bf k} + {\bf e}_l} \sum_{j=1}^3 \delta_{ij}\tilde{M}_{jl}^{{\bf k}} \bigg].
\end{align}
\end{subequations}
Then split~(\ref{T_2b}) into two terms, and simplify (\ref{T_2c}) and (\ref{T_2d}) to obtain
\begin{subequations}
\begin{align}
\label{T_3a}
u_i^{str} 
&\approx \frac{1}{3}\sum_{||{\bf k}||=0}^p 
\bigg[ (k_i + 1) \sum_{j=1}^3 (k_j + 1 + \delta_{ij})b^{{\bf k} + {\bf e}_i + {\bf e}_j} 
{\sum_{l=1}^3 (x_l - y_l)\tilde{M}_{jl}^{{\bf k}}} \\
\label{T_3b}
& - (k_i + 1) \sum_{j=1}^3 \sum_{l=1}^3 (k_j + 1 + \delta_{ij})b^{{\bf k} + {\bf e}_i + {\bf e}_j - {\bf e}_l}\tilde{M}_{jl}^{{\bf k}} \\
\label{T_3c}
& + (k_i + 1) \sum_{j=1}^3 \sum_{l=1}^3 \delta_{jl} b^{{\bf k} + {\bf e}_i + {\bf e}_j - {\bf e}_l}\tilde{M}_{jl}^{{\bf k}} \\
\label{T_3d}
& + \sum_{j=1}^3 (k_j + 1)b^{{\bf k} + {\bf e}_j} \tilde{M}_{ji}^{{\bf k}} + \sum_{l=1}^3 (k_l + 1)b^{{\bf k} + {\bf e}_l} \tilde{M}_{il}^{{\bf k}} \bigg].
\end{align}
\end{subequations}
Defining
$\tau_j = \sum\limits_{l=1}^3 (x_l - y_l)\tilde{M}_{jl}^{{\bf k}}$ in~(\ref{T_3a}),
simplifying~(\ref{T_3c}),
and 
combining the terms in~(\ref{T_3d}), we have
\begin{subequations}
\begin{align}
u_i^{str}
&\approx \frac{1}{3}\sum_{||{\bf k}||=0}^p 
\bigg[ (k_i + 1) \sum_{j=1}^3 (k_j + 1 + \delta_{ij})b^{{\bf k} + {\bf e}_i + {\bf e}_j} \tau_j \\
& - (k_i + 1) \sum_{j=1}^3 \sum_{l=1}^3 (k_j + 1 + \delta_{ij})b^{{\bf k} + {\bf e}_i + {\bf e}_j - {\bf e}_l}\tilde{M}_{jl}^{{\bf k}} \\
& + (k_i + 1) b^{{\bf k} + {\bf e}_i} { \sum_{j=1}^3 \tilde{M}_{jj}^{{\bf k}}}
+ \sum_{j=1}^3 (k_j + 1)b^{{\bf k} + {\bf e}_j} { (\tilde{M}_{ji}^{{\bf k}} + \tilde{M}_{ij}^{{\bf k}})} \bigg].
\end{align}
\end{subequations}
Next let
$m^{{\bf k}} = \sum\limits_{j=1}^3 \tilde{M}_{jj}^{{\bf k}},~
m_{ij}^{{\bf k}} = \tilde{M}_{ij}^{{\bf k}} + \tilde{M}_{ji}^{{\bf k}}$,
and
restore the arguments ${\bf x}^m, C$,
so that the stresslet particle-cluster approximation is
\begin{subequations}
\label{T_final}
\begin{align}
u_i^{str}({\bf x}^m,C)
&\approx \frac{1}{3}\sum_{||{\bf k}||=0}^p 
\bigg[ (k_i + 1) \sum_{j=1}^3 (k_j + 1 + \delta_{ij})b^{{\bf k} + {\bf e}_i + {\bf e}_j} \tau_j(C) \\
& - (k_i + 1) \sum_{j=1}^3 \sum_{l=1}^3 (k_j + 1 + \delta_{ij})b^{{\bf k} + {\bf e}_i + {\bf e}_j - {\bf e}_l}\tilde{M}_{jl}^{{\bf k}}(C) \\
& + (k_i + 1) b^{{\bf k} + {\bf e}_i} m^{{\bf k}}(C) + \sum_{j=1}^3 (k_j + 1)b^{{\bf k} + {\bf e}_j} m_{ij}^{{\bf k}}(C) \bigg].
\end{align}
\end{subequations}
In summary,
the particle-cluster approximation~(\ref{pc_interaction}) is computed 
using~(\ref{S_final}) for the Stokeslet part
and 
(\ref{T_final}) for the stresslet part,
where the coefficients $b^{\bf k}$ are computed using the recurrence relation~(\ref{recurrence_1}).
The introduction of the quantities $\sigma^{\bf k}, \tau^{\bf k}, m^{\bf k}, m_{ij}^{\bf k}$
enables a more efficient computation.
The operation count for the treecode is still $O(N\log N)$,
but we observe empirically that the CPU run time is reduced
since this procedure avoids explicitly forming the Taylor coefficient tensors
$a_{ij}^{\bf k}, \tilde{a}_{ijl}^{\bf k}$.

 
\section{Description of Treecode}

The treecode implementation starts by inputting the particle positions and weights,
and 
building a hierarchical tree of particle clusters~\cite{BarnesHut}.
The root cluster is a cube containing all the source particles. 
The root is bisected along the Cartesian axes and the eight children become subclusters of the root. 
The child clusters are similarly bisected and the process continues until a cluster contains fewer than $N_0$ particles, 
where $N_0$ is a user-specified parameter. 
Each cluster has a data structure containing necessary information, 
e.g. pointers to the particles belonging to the cluster, 
coordinates of the cluster center, 
pointers to the children of the cluster, 
and so on. 

The procedure is outlined in Algorithm \ref{treecode_algorithm}. 
The code uses a multipole acceptance criterion (MAC) to determine whether a 
given target particle ${\bf x}^m$
and 
source cluster $C$ are well-separated.
The criterion for being well-separated is 
\begin{equation}
\frac{r}{R} \leq \theta, \label{MAC}
\end{equation}
where as shown in Figure \ref{fig:particle_cluster},
$r$ is the cluster radius, 
$R$ is the distance between the particle and the cluster center ${\bf y}_c$,
and 
$\theta$ is a user-specified parameter which together with the order $p$ controls the 
approximation error.
The code cycles through the target particles,
and
each particle interacts with source clusters starting at the root.
If the MAC~(\ref{MAC}) is satisfied,
the particle-cluster approximation is computed as explained above.
If the MAC is not satisfied, the code checks the children of the cluster, 
or if the cluster is a leaf (no children), direct summation is performed,
again using the efficient formulation explained above.
This structure follows the original Barnes-Hut treecode algorithm~\cite{BarnesHut},
modified to accommodate higher-order particle-cluster approximations of the
Stokeslet and stresslet kernels.

\begin{algorithm}
\caption{treecode}\label{treecode_algorithm}
\begin{algorithmic}[1]
\State program {\bf main} 
\State \,\,\,\,\,\,\,\, input particle positions ${\bf x}^n$ and weights $f^n_j, h^n_j, \nu_j^n$,
treecode parameters $p, \theta, N_0$
\State \,\,\,\,\,\,\,\, build tree, compute cluster moments $M_j^{\bf k}(C), \tilde{M}_{jl}^{\bf k}(C)$
\State \,\,\,\,\,\,\,\, for $n = 1:N$
\State \,\,\,\,\,\,\,\,\,\,\,\,\,\,\,\,\, {\bf compute-velocity} ($\bf{ x}^n$, root-cluster)
\State subroutine {\bf compute-velocity} (${\bf x}$, $C$)
\State \,\,\,\,\,\,\,\, if MAC is satisfied
\State \,\,\,\,\,\,\,\,\,\,\,\, compute particle-cluster interaction by far-field Taylor approximation
\State \,\,\,\,\,\,\,\, else
\State \,\,\,\,\,\,\,\,\,\,\,\, if $C$ is a leaf, compute particle-cluster interaction by direct sum
\State \,\,\,\,\,\,\,\,\,\,\,\, else
\State \,\,\,\,\,\,\,\,\,\,\,\,\,\,\,\, for each child $C^\prime$ of $C$
\State  \,\,\,\,\,\,\,\,\,\,\,\,\,\,\,\,\,\,\,\,\,\,\,\, \,\,\, {\bf compute-velocity} ($\bf x$, $C^\prime$)
\end{algorithmic}
\end{algorithm}

The treecode algorithm was programmed in C++ 
and
compiled using the Intel icpc compiler with -O2 optimzation flag.
The source code is available for download~\cite{github}.
The computations were performed on the University of Wisconsin-Milwaukee 
Mortimer Faculty Research Cluster  which has 55 standard compute nodes
and
each node is a Dell PowerEdge R430 server with 
two 12-core Intel Xeon E5-2680 v3 processors at 2.50GHz and 64 GB RAM.
Serial computations were done on one core
and
parallel computations were done using MPI with each process running on one core.


\section{Numerical Results}

We present results for two test cases.
The first test case has Stokeslet and stresslet particles on the surface of a unit sphere,
as in boundary element simulations of exterior Stokes flow.
The particle distribution is given by triangulating the sphere as follows.
Starting from an icosahedron with 20 triangular faces,
the faces are refined by connecting the centers of the edges,
resulting in triangulations with $N = 20\,\cdot\,4^L$ faces for $L$ levels of refinement.
The particles are obtained by projecting the triangle centroids onto the sphere.
In this test case each cluster in the tree is shrunk to the bounding rectangular box 
containing its particles.
The second test case follows~\cite{Klinteberg2017}
and
considers Stokeslet particles randomly distributed in cubes of different sizes,
with a fixed particle number density.
In both test cases the particle weights $f^n_j$, $ h^n_j$ are random numbers in $[-1, 1]$. 

There are three user-specified parameters required for the treecode. 
The parameter $N_0$ is the maximum number of particles in a leaf of the tree;
the value $N_0 = 2000$ is used throughout this work.
The other parameters are the
order of approximation $p$
and
MAC parameter $\theta$;
we will vary these parameters to investigate their effect on the code's performance.
The relative error in velocity~(\ref{StoStre}) is defined by
\begin{equation}
E = 
\left(\frac{\displaystyle \sum_{n = 1}^{N}|{\bf u}^d({\bf x}^n) - {\bf u}^t({\bf x}^n)|^2}
{\displaystyle \sum_{n = 1}^{N}|{\bf u}^d({\bf x}^n)|^2}\right)^{\!\!1/2}, 
\label{L2Err}
\end{equation}
where ${\bf u}^d$ is the velocity obtained by direct summation,
and 
${\bf u}^t$ is the treecode approximation.
The CPU time is given in units of seconds (s).
The following subsection presents serial results for each test case,
followed by a subsection with parallel results for the second test case.


\subsection{Serial computations}

\subsubsection{Test Case 1: particles on the surface of a sphere}

\begin{figure}[htb]
\centering
\includegraphics[width=0.8\linewidth]{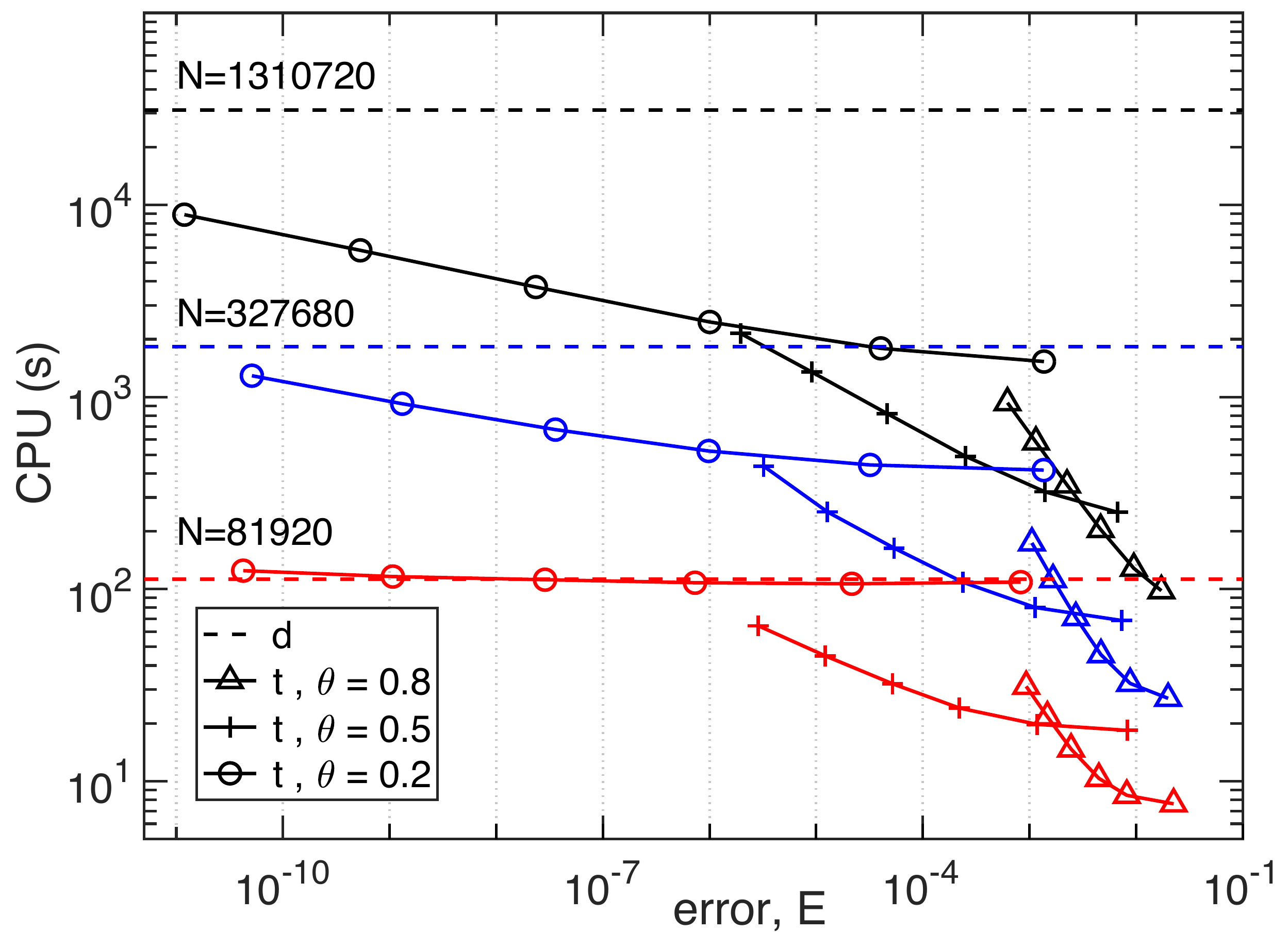}
\caption{Test case 1,
Stokeslets and stresslets on a sphere, 
CPU time versus error $E$, 
system size 
$N = 1310720$ (upper, black),
$N = 327680$ (middle, blue),
$N = 81920$ (lower,red),  
direct sum (d, dashed),
treecode (t, solid),
MAC parameter $\theta = 0.8\,(\triangle), 0.5\,(+), 0.2\,(\bigcirc)$ (decreasing right to left),
order $p = 0:2:10$ (increasing right to left).}
\label{fig:TreecodePerformance_SCATTER125}
\end{figure}

The first test case has Stokeslet and stresslet particles on a unit sphere
at locations determined by the icosahedral triangulation as explained above.
Figure~\ref{fig:TreecodePerformance_SCATTER125} 
plots the treecode CPU time versus error $E$ for systems of size 
$N = 82920, 317680, 1310720$,
with Taylor approximation order $p = 0:2:10$ (increasing from right to left),
and
MAC parameter $\theta = 0.8,0.5,0.2$ (decreasing from right to left).
As expected,
smaller error $E$ is attained by increasing the order $p$ and decreasing the MAC parameter $\theta$,
but this increases the CPU time.
Note that as $p$ increases, the error decreases more rapidly for smaller $\theta$.
Also note that for a given order $p$ and MAC parameter $\theta$,
the error is relatively insensitive to the system size $N$.
We can distinguish three regimes;
for low accuracy $\theta = 0.8$ is most efficient,
for medium accuracy $\theta = 0.5$ is most efficient,
and
for high accuracy $\theta = 0.2$ is most efficient.
Figure \ref{fig:TreecodePerformance_SCATTER125} also shows the direct sum CPU time;
as the system size increases,
the treecode becomes more efficient in comparison with direct summation.

The dependence of the error and CPU time on the MAC parameter $\theta$
can be explained as follows.
Choosing a smaller $\theta$ has two effects,
(i) the code descends deeper into the tree,
so the CPU time increases,
but the clusters have smaller radius,
so the Taylor approximation is more accurate, 
and 
(ii) there is higher likelihood the code will reach the leaf clusters of the tree,
which again increases the CPU time,
but in that case direct summation is performed,
which incurs no error.

Table \ref{tab:TestCase1ScatterData} gives the speedup
(ratio $d/t$ of CPU times for direct sum and treecode),
and the error $E$.
The treecode achieves higher speedup as the system size increases;
for example,
a medium accuracy $E \approx {\rm 5e{-}05}$ is attained using 
MAC parameter $\theta = 0.5$ and order $p=6$;
in this case the treecode is 
3.5 times faster than direct summation for $N=81920$,
11 times faster for $N=327680$
and
38 times faster for $N=1310720$.

\begin{table}[htb]
\centering
\begin{tabular}{|c||r|r|c||r|r|c||r|r|c|}
\hline
\multicolumn{1}|{c||}{} &
\multicolumn{3}{c||}{(a) $N = 81920$} & 
\multicolumn{3}{c||}{(b) $N = 327680$} & 
\multicolumn{3}{c|}{(c) $N = 1310720$} \\
\hline
$\theta$ & $p$ & $d/t$ & $E$ & $p$ & $d/t$ & $E$ & $p$ & $d/t$ & $E$ \\
\hline
      & 0   & 1.03 & 8.3e-04 & 0   & 4.38 & 1.4e-03 & 0   & 20.43 & 1.4e-03 \\
      & 2   & 1.05 & 2.2e-05 & 2   & 4.11 & 3.2e-05 & 2   & 17.44 & 4.0e-05 \\
      & 4   & 1.04 & 7.3e-07 & 4   & 3.48 & 9.8e-07 & 4   & 12.71 & 1.0e-06 \\
0.2 & 6   & 1.00 & 2.9e-08 & 6   & 2.69 & 3.6e-08 & 6   & 8.36   & 2.4e-08 \\
      & 8   & 0.97 & 1.1e-09 & 8   & 1.97 & 1.3e-09 & 8   & 5.38   & 5.3e-10 \\
      & 10 & 0.90 & 4.2e-11 & 10 & 1.41 & 5.1e-11 & 10 &  3.51  & 1.2e-11 \\
\hline
      & 0   & 6.12 & 8.9e-03 & 0   & 26.59 & 7.4e-03 & 0   & 124.49 & 6.7e-03 \\
      & 2   & 5.69 & 1.2e-03 & 2   & 22.75 & 1.1e-03 & 2   & 97.36   & 1.4e-03 \\
      & 4   & 4.68 & 2.2e-04 & 4   & 16.79 & 2.4e-04 & 4   & 63.67   & 2.5e-04 \\
0.5 & 6   & 3.50 & 5.2e-05 & 6   & 11.15 & 5.4e-05 & 6   & 38.12   & 4.6e-05 \\
      & 8   & 2.50 & 1.2e-05 & 8   & 7.23   & 1.3e-05 & 8   & 23.18   & 9.1e-06 \\
      & 10 & 1.74 & 2.9e-06 & 10 & 4.18   & 3.2e-06 & 10 & 14.59   & 2.0e-06 \\
\hline
      & 0   & 14.75 & 2.3e-02 & 0   & 67.63 & 2.0e-02 & 0   & 318.73 & 1.7e-02 \\
      & 2   & 13.32 & 8.2e-03 & 2   & 56.52 & 8.8e-03 & 2   & 243.78 & 9.5e-03 \\
      & 4   & 10.87 & 4.5e-03 & 4   & 25.66 & 4.7e-03 & 4   & 153.01 & 4.7e-03 \\
0.8 & 6   & 7.65   & 2.5e-03 & 6   & 25.66 & 2.7e-03 & 6   & 89.71   & 2.3e-03 \\
      & 8   & 5.19   & 1.5e-03 & 8   & 16.29 &1.7e-03 &  8   & 53.60   & 1.2e-03 \\
      & 10 & 3.61   & 9.3e-04 & 10 & 10.49 & 1.1e-03 & 10 & 33.48   & 6.2e-04 \\
\hline
\end{tabular}
\caption{Test case 1,
Stokeslets and stresslets on a sphere, 
treecode MAC parameter $\theta$,
order $p$,
speedup $d/t$ = ratio of direct sum and treecode CPU times,
error $E$,
(a) $N=81920$, $d$ = 112 s,
(b) $N = 327680$, $d$ = 1824 s,
(c) $N = 1310720$, $d$ = 31267 s.}
\label{tab:TestCase1ScatterData}
\end{table}

\begin{figure}[h!]
\centering
\includegraphics[width=0.8\linewidth]{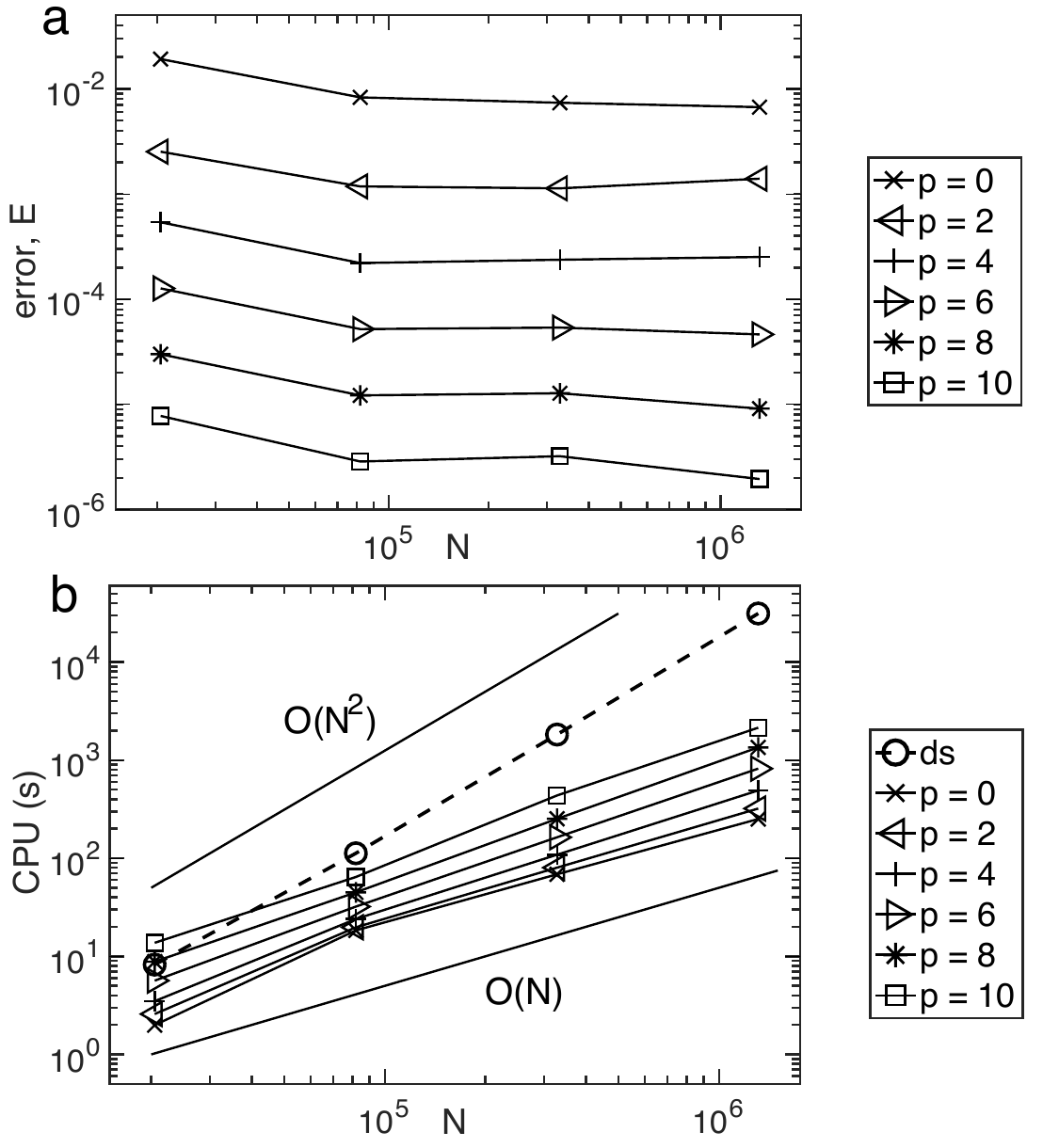}
\caption{Test case 1,
Stokeslets and stresslets on a sphere,
treecode MAC parameter $\theta = 0.5$,
order $p = 0:2:10$,
system size $N = 20480:1310720$, 
(a) error, $E$, 
(b) CPU time for direct sum (dashed), treecode (solid).}
\label{fig:TwoSphere_ErrorCPU}
\end{figure}

Figure \ref{fig:TwoSphere_ErrorCPU} shows results graphically 
for MAC parameter $\theta = 0.5$,
order $p=0:2:10$,
and
system size between $N=20480$ and $N=1310720$.
Figure \ref{fig:TwoSphere_ErrorCPU}a plots the treecode error $E$ versus $N$.
The error is relatively insensitive to the system size,
and for a given $N$, the error decreases with increasing order $p$. 
Figure \ref{fig:TwoSphere_ErrorCPU}b plots the CPU time versus $N$. 
The direct sum CPU time scales like $O(N^2)$,
while the treecode CPU time is consistent with $O(N \log N)$,
and
hence the treecode is faster than direct summation except for 
small system size $N$ and large order $p$. 
In the remainder of this section we use MAC parameter $\theta = 0.5$.
   
Table~\ref{tab:TreecodePerformance_Memory}
displays the peak memory used by the treecode as a function of system size $N$ for order $p = 0 : 2 : 10$.
The memory used by direct summation is also given.  
The memory usage statistics were obtained using the Valgrind massif analysis tool
(valgrind.org).
The treecode and direct summation both store the particles in arrays of size $O(N)$,
but the treecode uses additional memory of size $O((N/N_0)p^3)$ for the cluster moments,
where the factor $N/N_0$ represents the number of clusters in the tree
and
the factor $p^3$ is the memory associated with the moments of a cluster.
Recall that the number of particles in a leaf cluster is set here at $N_0 = 2000$.

In the current treecode implementation,
the additional memory used for the moments was reduced as follows.
Note that for a given cluster $C$ and indices $j,l = 1:3$,
the moments $M_j^{\bf k}(C), \tilde{M}_{jl}^{\bf k}(C)$ could be stored in
square three-dimensional arrays of size $p^3$,
corresponding to the index ${\bf k} = (k_1,k_2,k_3)$.
However since we only need the indices with $||{\bf k}|| = k_1 + k_2 + k_3 = 0:p$,
a large portion of these arrays would be empty.
So instead the moments are stored in one-dimensional arrays with no empty space,
by accessing the indices $(k_1,k_2,k_3)$ in a fixed order;
we refer to this as flattening the moment arrays.
While the treecode memory usage still has a term scaling like $O((N/N_0)p^3)$,
the prefactor is greatly reduced.
Table~\ref{tab:TreecodePerformance_Memory}
shows that in this range of system size $N$
and
order $p \le 10$,
except for one case with the smallest $N$
and
largest $p$,
the treecode uses less than twice as much memory as direct summation.
 
\begin{table}[htb]
\centering
\begin{tabular}{|c|c|c|c|c|}
\hline
 & $N = 20480$ & $N = 81920$ & $N = 327680$ & $N = 1310720$ \\
\hline
direct sum & 3.3 & 13.1 & 52.4 & 209.8 \\
\hline
treecode order $p$ & & & & \\
\hline
  0 &   3.3 & 13.1 & 52.6 & 210.7 \\
  2 &   3.4 & 13.2 & 53.0 & 214.3 \\
  4 &   3.7 & 13.6 & 54.2 & 224.1 \\
  6 &   4.3 & 14.2 & 56.6 & 243.4 \\
  8 &   5.4 & 15.2 & 60.6 & 275.2 \\
10 &   7.6 & 16.6 & 66.6 & 322.8 \\
\hline
\end{tabular}
\caption{
Test case 1,
Stokeslets and stresslets on a sphere,
memory usage (MB) for system size $N$,
direct sum, treecode with order $p$.}
\label{tab:TreecodePerformance_Memory}
\end{table}

\subsubsection{Test Case 2: particles in a cube}

The second test case has Stokeslet particles located randomly in a cube of side length $L$,
with number density $N/L^3 = 2500$~\cite{Klinteberg2017}.
Table~\ref{tab:TestCase2ScatterData} 
shows the speedup (ratio $d/t$ of CPU times for direct sum and treecode),
and the error $E$,
for system size $N=125K$ and $N=1000K$.
The trends in error and CPU time with respect to 
MAC parameter $\theta$ and order $p$ are similar to test case 1.
The speedup in this case is somewhat less than in the previous case,
because here the system sizes are smaller
and
the stresslet part of the sum is omitted;
nonetheless the treecode is faster than direct summation except for one case with
small $\theta$ and large $p$.

\begin{table}[htb]
\centering
\begin{tabular}{|c||r|r|c||r|r|c|}
\hline
\multicolumn{1}{|c||}{} &
\multicolumn{3}{c||}{(a) $N = 125$K} &
\multicolumn{3}{c|}{(b) $N = 1000$K} \\
\hline
$\theta$ & $p$ & $d/t$ & $E$ & $p$ & $d/t$ & $E$ \\
\hline
      & 0   & 1.14 & 1.1e-02 & 0   & 3.46 & 4.3e-02 \\
      & 2   & 1.13 & 2.2e-04 & 2   & 3.37 & 7.1e-04 \\
0.2 & 4   & 1.11 & 4.4e-06 & 4   & 3.18 & 1.5e-05  \\
      & 6   & 1.08 & 1.1e-07 & 6   & 2.89 & 3.9e-07 \\
      & 8   & 1.02 & 3.2e-09 & 8   & 2.50 & 1.1e-08 \\
      & 10 & 0.95 & 9.6e-11 & 10 & 2.08 & 3.2e-10 \\
\hline
      & 0   & 5.50 & 1.1e-01 & 0   & 33.89 & 1.5e-01 \\
      & 2   & 5.28 & 8.2e-03 & 2   & 31.39 & 1.6e-02 \\
0.5 & 4   & 4.79 & 1.1e-03 & 4   & 26.86 & 2.0e-03 \\
      & 6   & 4.17 & 1.7e-04 & 6   & 21.13 & 3.3e-04 \\
      & 8   & 3.40 & 3.0e-05 & 8   & 13.77 & 5.7e-05 \\
      & 10 & 2.64 & 5.5e-06 & 10 & 11.28 & 1.0e-05 \\
\hline
      & 0   & 16.94 & 2.1e-01 & 0   & 123.78 & 2.7e-01 \\
      & 2   & 15.85 & 4.3e-02 & 2   & 111.61 & 7.0e-02 \\
0.8 & 4   & 13.81 & 1.4e-02 & 4   & 90.79   & 2.3e-02 \\
      & 6   & 11.32 & 5.4e-03 & 6   & 67.34   & 9.0e-03 \\
      & 8   & 8.77   & 2.4e-03 & 8   & 47.13   & 4.1e-03 \\
      & 10 & 6.40   & 1.1e-03 & 10 & 32.54   & 1.9e-03 \\
\hline
\end{tabular}
\caption{Test case 2,
random Stokeslets in a cube with number density $N/L^3 = 2500$~\cite{Klinteberg2017},
treecode MAC parameter $\theta$,
order $p$,
speedup $d/t$ = ratio of direct sum and treecode CPU times,
error $E$,
(a) $N = 125$K, $d$ = 150 s,
(b) $N = 1000$K, $d$ = 9873 s.}
\label{tab:TestCase2ScatterData}
\end{table}

Table~\ref{tab:TreecodePerformance_Memory_2} presents the memory usage in test case 2.
As before, the treecode uses less than twice as much memory as direct summation.
The relatively low memory usage of the treecode is an advantage in parallel simulations,
where it enables a simple replicated data approach as shown below.

\begin{table}[htb]
\centering
\begin{tabular}{|c|c|c|}
\hline
 & $N = 125K$ & $N = 1000K$ \\
\hline
direct sum & 11.0 & 88.0 \\
\hline
treecode order $p$ & & \\
\hline
  0 & 14.0 & 112.4  \\
  2 & 14.1 & 112.8  \\
  4 & 14.3 & 114.0 \\
  6 & 14.6 & 116.4  \\
  8 & 15.0 & 120.4 \\
10 & 15.8 & 126.3 \\
\hline
\end{tabular}
\caption{
Test case 2,
random Stokeslets in a cube with number density $N/L^3 = 2500$~~\cite{Klinteberg2017},
memory usage (MB) for system size $N$,
direct sum, treecode with order $p$.}
\label{tab:TreecodePerformance_Memory_2}
\end{table}

 
\subsection{Parallel computations}

The parallel treecode implementation used here relies on the observation
that the target particle computations in Algorithm~\ref{treecode_algorithm} are independent of each other,
enabling a replicated data approach~\cite{Geng2013}.
The procedure is sketched in Algorithm~\ref{parallel_treecode}
and
was implemented using MPI.
The idea is that the particle array of length $N$ is divided into $np$ segments of length $N/np$, 
where $np$ is the chosen number of processes,
and
each segment is assigned to one process
(recall that each process is running on a single core).
Each process receives a copy of the entire particle array,
builds a local copy of the tree,
and
computes the cluster moments.
The processes run concurrently
and
each one computes the induced velocity at its assigned target particles. 
The scheme assumes that the entire particle array and tree structure fit into the memory of each core,
which is facilitated by the treecode's relatively low memory usage;
this is not an issue for the examples considered here,
but in case of a much larger system size where this assumption doesn't hold, 
a distributed memory approach would be required.

The parallel treecode performance is demonstrated below for test case 2.
First we consider strong scaling (fixed $N$),
and
then weak scaling (increasing $N$),
as the number of processes $np$ increases.
We also implemented a parallel direct sum using the same replicated data approach
for comparison with the parallel treecode.
  
\begin{algorithm}
\caption{parallel treecode}\label{parallel_treecode}
\begin{algorithmic}[1]
\State in main process
\State \,\,\,\,\,\,\,\,\,\, input particle positions ${\bf x}^n$, weights $f_j^n, h_j^n$, 
treecode parameters $p, \theta, N_0$
\State \,\,\,\,\,\,\,\,\,\, broadcast particle array to each process
\State in each process
\State \,\,\,\,\,\,\,\,\,\, build local copy of tree, compute cluster moments
\State \,\,\,\,\,\,\,\,\,\, use treecode to compute induced velocity at assigned target particles
\State \,\,\,\,\,\,\,\,\,\, send result to main process
\end{algorithmic}
\end{algorithm}

\subsubsection{Strong scaling}

Table \ref{table:parallel} shows results for test case 2
with $N = 1000$K random Stokeslets in a cube.
The treecode parameters are $\theta = 0.5, p = 6$, 
yielding error $E = 3.1{\rm e{-}04}$.
The Table displays the 
CPU time,
ratio of CPU time for 1~process and $np$ processes,
and
parallel efficiency,
for the direct sum ($d$) and treecode ($t$),
and 
finally the speedup due to using the treecode $(d/t)$, 
up to $np = 32$ processes. 
The parallel direct sum reduces the CPU time from 
9707~s on 1 process to
350.7~s on 32 processes,
for parallel efficiency 86.5\%.
As expected the treecode CPU times are smaller,
and
the parallel performance is almost as good;
the CPU time is reduced from
336.6~s on 1 process to
13.6~s on 32 processes,
for parallel efficiency 77.3\%.
The treecode is 28 times faster than direct summation on 1~process
and
25 times faster on 32 processes.
  
\begin {table}[htb]  
\begin{center}
\begin{tabular}{|r || r | r | c || r | r | c || r|}
\hline
$np$ & $d$ CPU (s) & $d_1/d_{np}$ & $d$ PE $(\%)$ & $t$ CPU (s) & $t_1/t_{np}$ & $t$ PE $(\%)$ & $d/t$ \\ 
\hline
1   & 9707.0 & 1.00   & 100.0 & 336.6 &  1.00 & 100.0 & 28.8 \\ \hline
2   & 4911.4 & 1.98   & 98.8   & 176.3 &  1.91 &   95.4 & 27.9 \\ \hline
4   & 2451.2 & 3.96   & 99.0   & 86.6   &  3.89 &   97.2    & 28.3 \\ \hline
8   & 1345.6 & 7.21   & 90.2   &  48.1  &  7.00 &   87.6    & 28.0 \\ \hline
16 & 702.0   & 13.83 & 86.4   & 25.6   & 13.16 &  82.2    & 27.4 \\ \hline
32 &  350.7  & 27.68 & 86.5   & 13.6   & 24.75 &  77.3    & 25.8 \\ \hline
\end{tabular}
\caption{Test case 2,
random Stokeslets in a cube,
parallel strong scaling,
system size $N=1000$K,
treecode parameters $\theta = 0.5$, $p = 6$,
error $E = 3.1{\rm e{-}04}$,
number of MPI processes ($np$),
CPU time ($d$ = direct sum, $t$ = treecode),
ratio of CPU time for one process and $np$ processes ($d_1/d_{np}, t_1/t_{np}$),
parallel efficiency (PE, ratio/$np$),
treecode speedup ($d/t$).}
\label{table:parallel}
\end{center}
\end{table}

\subsubsection{Weak scaling}

Table \ref{table:parallel_weakscale} shows results for test case 2
starting with $N = 125$K particles on 1 process, 
then doubling the number of particles 
and 
processes until reaching $N = 4000$K particles on 32 processes.
The box size increases so that the number density is $N/L^3 = 2500$~\cite{Klinteberg2017}.
The treecode parameters are $\theta = 0.5, p=6$,
yielding error $E \le 3.7{\rm e-04}$.
The results show that with each doubling of $N$ and $np$,
the direct sum CPU time approximately doubles,
while the treecode CPU time increases more slowly.
Hence the treecode performance improves as the system size increases;
with $N=125$K on 1 process, 
the treecode is 6 times faster than direct summation,
but with $N=4000$K on 32 processes,
the treecode is 82 times faster.  
       
\begin {table}[htb]  
\begin{center}
\begin{tabular}{|r | r | r | r | r|}
\hline
$N$ & $np$ & $d$ CPU (s) & $t$ CPU (s) & $d/t$ \\ \hline
   125K &   1 &  145.5  & 23.7  & 6.1  \\ \hline
   250K &   2 &  291.7  & 27.2  & 10.7 \\ \hline
   500K &   4 &   610.3 & 32.9  & 18.6 \\ \hline
 1000K &   8 & 1345.7 & 48.1  & 28.0 \\ \hline
 2000K & 16 & 2804.9 & 56.7  & 49.5 \\ \hline 
 4000K & 32 & 5565.9 & 67.7  & 82.3 \\ \hline  
\end{tabular}
\caption{Test case 2,
random Stokeslets in a cube,
parallel weak scaling, 
system size ($N$),
number of MPI processes ($np$),
treecode parameters $\theta = 0.5$, $p = 6$, 
error $E \le 3.7{\rm e{-}04}$,
CPU time ($d$ = direct sum, $t$ = treecode),
treecode speedup $d/t$.}
\label{table:parallel_weakscale}
\end{center}
\end{table}
    
  
\section{Summary}

We presented a treecode algorithm for computing the velocity induced by a collection of 
Stokeslets and stresslets in 3D flow.
The method uses a far-field Cartesian Taylor approximation to compute well-separated particle-cluster interactions.
Expressions were derived for the Taylor coefficients of the Stokeslet and stresslet kernels
in terms of the Taylor coefficients of the Coulomb potential,
and
these expressions enable an efficient computation of higher-order approximations.
Numerical results were presented for 
icosahedral particles on the surface of a sphere,
and 
random particles in a cube~\cite{Klinteberg2017}.
For a given level of accuracy, 
the treecode CPU time scales like $O(N \log N )$, 
where $N$ is the number of particles, 
and a substantial speedup over direct summation is achieved for large systems. 
The memory usage increases with the system size $N$ and Taylor approximation order $p$,
but for the range of parameters considered here,
in most cases the treecode used less than twice as much memory as direct summation.
A relatively straightforward parallel treecode implementation was demonstrated. 

It is beyond the scope of the present work to make a detailed performance comparison
with other methods for fast summation of Stokeslets and stresslets
such as
the Fast Multipole Method (FMM)~\cite{Gimbutas2015,MalhotraBiros2015,MalhotraBiros2016}
and
the Spectral Ewald (SE) method~\cite{Klinteberg2017}.
These methods have demonstrated excellent performance in terms of accuracy and efficiency.
Yet the treecode may be an attractive option in some cases due to its
relatively simple algorithmic structure and low memory usage,
which together can enhance parallel efficiency.

Our simulations used representative values of the treecode parameters
(order $p$, MAC parameter $\theta$, maximum number of particles in a leaf $N_0$)
and
one future goal is to gain efficiency by tuning their values.
There are several other directions for future work.
The present approach can be extended to treat regularized Stokeslets and stresslets,
which may help accelerate biofluid applications using those 
kernels~\cite{Cortez2001,cortez-fauci-medovikov-05,Rostami2016}.
Another goal is to apply the treecode in
boundary element simulations of 
Stokes-Darcy porous medium flow~\cite{tlupova-cortez-09,boubendir-tlupova-13}
and
Stokes flow around solid bodies~\cite{Pozrikidis1992}.

\section*{Acknowledgments}

We thank the reviewers for suggestions that helped improve the article,
and
Leighton Wilson for suggesting the method used in reducing the memory associated with particle clusters. 



\begin{thebibliography}{99}

\bibitem{Klinteberg2016}
L. af Klinteberg and A.-K. Tornberg,
{\it A fast integral equation method for solid particles in viscous flow using quadrature by expansion},
J. Comput. Phys., 326 (2016), pp.~420--445.

\bibitem{Klinteberg2017}
L. af Klinteberg, D.S. Shamshirgar and A.-K. Tornberg,
{\it Fast Ewald summation for free-space Stokes potentials},
Res. Math. Sci., 4 (2017)4:1, DOI 10.1186/s40687-016-0092-7.

\bibitem{ambrose-siegel-tlupova-13}
D.M. Ambrose, M. Siegel and S. Tlupova,
{\it A small-scale decomposition for 3D boundary integral computations with surface tension},
J. Comput. Phys., 247 (2013), pp.~168-191.

\bibitem{BarnesHut}
J.E. Barnes and P. Hut,
{\it A hierarchical $O(N\log N)$ force-calculation algorithm},
Nature, 324 (1986), pp.~446--449.

\bibitem{BeatusTlustyBarZiv}
T. Beatus, T. Tlusty and R. Bar-Ziv,
{\it Phonons in a one-dimensional microfluidic crystal},
Nature Physics, 2 (2006), pp.~743--748.



\bibitem{boubendir-tlupova-13}
Y. Boubendir and S. Tlupova,
{\it Domain decomposition methods for solving Stokes-Darcy problems with boundary integrals},
SIAM J. Sci. Comput., 35 (2013), pp.~B82-B106.

\bibitem{CoronaGreengardRachhVeerapaneni2017}
E. Corona, L. Greengard, M. Rachh and S. Veerapaneni,
{\it An integral equation formulation for rigid bodies in Stokes flow in three dimensions},
J. Comput. Phys., 332 (2017), pp.~504--519.

\bibitem{Cortez2001}
R. Cortez,
{\it The method of regularized Stokeslets},
SIAM J. Sci. Comput., 23 (2001), pp.~1204--1225.

\bibitem{cortez-fauci-medovikov-05}
R. Cortez, L. Fauci and A. Medovikov,
{\it The method of regularized Stokeslets in three dimensions:
Analysis, validation, and application to helical swimming},
Phys. Fluids, 17 (2005), 031504.

\bibitem{Draghicescu1995}
C.I. Draghicescu and M. Draghicescu,
{\it A fast algorithm for vortex blob interactions},
J. Comput. Phys., 116 (1995), pp.~69--78.

\bibitem{Drescher2010}
K. Drescher, R.E. Goldstein, N. Michel, M. Polin and I. Tuval,
{\it Direct measurement of the flow field around swimming microorganisms},
Phys. Rev. Lett., 105 (2010), 168101.

\bibitem{Duan2001}
Z.-H. Duan and R. Krasny,
{\it An adaptive treecode for computing nonbonded potential energy in classical molecular systems},
J. Comput. Chem., 22 (2001), pp.~184--195.


\bibitem{Darden1995}
U. Essmann, L. Perera, M.L. Berkowitz, T. Darden, H. Lee and L.G. Pedersen,
{\it A smooth particle mesh Ewald method},
J. Chem. Phys., 103 (1995), pp.~8577--8593.


\bibitem{FMMLIB3D}
FMMLIB3D 1.2,
www.cims.nyu.edu/cmcl/fmm3dlib/fmm3dlib.html (2012)

\bibitem{Fu2000} 
Y. Fu and G.J. Rodin,
{\it Fast solution methods for three-dimensional Stokesian many-particle problems},
Commun. Numer. Meth. Engng., 16 (2000), pp.~145--149.

\bibitem{Geng2013}
W.-H. Geng and R. Krasny,
{\it A treecode-accelerated boundary integral Poisson-Boltzmann solver for solvated biomolecules},
J. Comput. Phys., 247 (2013), pp.~62--78.

\bibitem{Gimbutas2015} 
Z. Gimbutas and L. Greengard,
{\it Simple FMM libraries for electrostatics, slow viscous flow, and frequency-domain wave propagation},
Commun. Comput. Phys., 18 (2015), pp.~516--528.

\bibitem{github}
github.com/Treecodes/stokes-treecode

\bibitem{Greengard1987}
L. Greengard and V. Rokhlin,
{\it A fast algorithm for particle simulations},
J. Comput. Phys., 73 (1987), pp.~325--348.

\bibitem{Hockney1988}
R.W. Hockney and J.H. Eastwood,
Computer Simulation Using Particles,
Taylor \& Francis, 1988.


\bibitem{Wang2011}  
R. Krasny and L. Wang,
{\it Fast evaluation of multiquadric RBF sums by a Cartesian treecode},
SIAM J. Sci. Comput., 33 (2011), pp.~2341--2355.

\bibitem{LiJohnstonKrasny2009}  
P. Li, H. Johnston and R. Krasny,
{\it A Cartesian treecode for screened Coulomb interactions},
J. Comput. Phys., 228 (2009), pp.~3858--3868.

\bibitem{Lindsay2001} 
K. Lindsay and R. Krasny,
{\it A particle method and adaptive treecode for vortex sheet motion in three-dimensional flow},
J. Comput. Phys., 172 (2001), pp.~879--907.

\bibitem{Liu2008} 
Y.J. Liu,
{\it A new fast multipole boundary element method for solving 2D Stokes flow problems 
based on a dual BIE formulation},
Eng. Anal. Bound. Elem., 32 (2008), pp.~139--151.

 
\bibitem{MalhotraBiros2015}
D. Malhotra and G. Biros,
{\it PVFMM: A parallel kernel independent FMM for particle and volume potentials},
Commun. Comput. Phys., 18 (2015), pp.~808--830.

\bibitem{MalhotraBiros2016}
D. Malhotra and G. Biros,
{\it Algorithm 967: A distributed-memory Fast Multipole Method for volume potentials},
ACM Trans. Math. Soft., 43(2) (2016), 17.

 \bibitem{NitscheParthasarathi}
L.C. Nitsche and P. Parthasarathi,
{\it Cubically regularized Stokeslets for fast particle simulations of low-Reynolds-number drop flows},
Chem. Eng. Comm., 197 (2010), pp.~18--38.
 
\bibitem{PignatelNicolasGuazzelliSaintillan}
F. Pignatel, M. Nicolas, E. Guazzelli and D. Saintillan,
{\it Falling jets of particles in viscous fluids},
Phys. Fluids, 21 (2009), 123303.

\bibitem{Pozrikidis1992}  
C. Pozrikidis,
Boundary Integral and Singularity Methods for Linearized Viscous Flow,
Cambridge University Press, 1992.


\bibitem{Rostami2016}  
M.W. Rostami and S.D. Olson, 
{\it Kernel-independent fast multipole method within the framework of regularized Stokeslets},
J. Fluid. Struct., 67 (2016), pp.~60--84.

\bibitem{Saintillan2005}  
D. Saintillan, E., Darve and E. Shaqfeh,
{\it A smooth particle-mesh Ewald algorithm for Stokes suspension simulations: The sedimentation of fibers},
Phys. Fluids, 17 (2005), 033301/21.

\bibitem{Sierou2001} 
A. Sierou and J.F. Brady,
{\it Accelerated Stokesian Dynamics simulations},
J. Fluid Mech., 448 (2001), pp.~115--146.

\bibitem{Smith2009}
D.J. Smith,
{\it A boundary element regularized Stokeslet method applied to cilia- and flagella-driven flow},
Proc. R. Soc. A, 465 (2009), pp.~3605--3626.

\bibitem{tlupova-cortez-09}
S. Tlupova and R. Cortez,
{\it Boundary integral solutions of coupled Stokes and Darcy flows},
J. Comput. Phys., 228 (2009), pp.~158-179.

\bibitem{Tornberg2008} 
A.-K. Tornberg and L. Greengard,
{\it A fast multipole method for the three-dimensional Stokes equations},
J. Comput. Phys., 227 (2008), pp.~1613--1619.

\bibitem{VeerapaneniRahimianBirosZorin}
S.K. Veerapaneni, A. Rahimian, G. Biros and D. Zorin, 
{\it A fast algorithm for simulating vesicle flows in three dimensions},
J. Comput. Phys., 230 (2011), pp.~5610--5634.

\bibitem{Wang2007} 
H. Wang, T. Lei, J. Li, J. Huang and Z. Yao,
{\it A parallel fast multipole accelerated integral equation scheme for 3D Stokes equations},
Int. J. Numer. Meth. Engng., 70 (2007), pp.~812--839.

\bibitem{Wang2006} 
X. Wang, J. Kanpka, W. Ye, N.R. Aluru and J. White,
{\it Algorithms in Fast Stokes and its application to micromachined device simulation},
IEEE Trans. Comput. Aided Des. Integra. Circ. Syst., 25 (2006), pp.~248--257.

\bibitem{Ying2004}  
L. Ying, G. Biros and D. Zorin,
{\it A kernel independent adaptive fast multipole algorithm in two and three-dimensions},
J. Comput. Phys., 196 (2004), pp.~591--626.

%
%
%
%

\end{thebibliography}
\end{document}